\documentclass{article}
\usepackage{romp34}
\usepackage{amssymb,amscd}
\usepackage{amsmath}
\usepackage{bbm}
\newtheorem{proposition}{Proposition}

\input{prepictex.tex}
\input{pictex.tex}
\input{postpictex.tex}

\newcommand{\beq}{\begin{equation}}
\newcommand{\eeq}{\end{equation}}
\newcommand{\bpr}{\begin{proposition}}
\newcommand{\epr}{\end{proposition}}

\def\epf{\hfill$\square$}
\def\Box{\square}

\def\Z{\mathbb Z}
\def\R{\mathbb R}
\def\C{\mathbb C} 
 
\def\N{\mathbb N}

\def\O{{\mathcal O}}
\def\cO{\O}

\def\Prim{\operatorname{Prim}}
\title{Graph \boldmath$C^*$-algebras and\\ \boldmath$\Z_2$-quotients of quantum spheres} 
\author{
Piotr M.~Hajac\\
Mathematisches Institut, Universit\"at M\"unchen\\
Theresienstr.\ 39, M\"unchen, 80333, Germany\\
and\\
Instytut Matematyczny, Polska Akademia Nauk\\
ul.\ \'Sniadeckich 8, Warszawa, 00-950 Poland\\
and\\
Katedra Metod Matematycznych Fizyki, Uniwersytet Warszawski\\
ul.\ Ho\.za 74, Warszawa, 00-682 Poland \vspace{0mm}\\
http://www.fuw.edu.pl/$\!\widetilde{\phantom{m}}\!$pmh
\vspace*{2mm}\and
Rainer Matthes\thanks{~Speaker}\\
Fachbereich 2, TU Clausthal\\
Leibnizstr. 4, D-38678 Clausthal-Zellerfeld, Germany\\
e-mail: ptrm@pt.tu-clausthal.de
\vspace*{2mm}\and
Wojciech Szyma\'nski\\
School of Mathematical and Physical Sciences, University of Newcastle\\
Callaghan, NSW 2308, Australia\\
e-mail: wojciech@frey.newcastle.edu.au
}

\begin{document}
\maketitle
\begin{abstract}
We consider two $\Z_2$-actions on the Podle\'s generic 
quantum spheres. They
yield, as noncommutative quotient spaces, the Klimek-Lesniewski $q$-disc and the 
quantum real projective
space, respectively. The $C^*$-algebras of all 
these quantum spaces
are described as graph $C^*$-algebras. 
The $K$-groups of the thus presented $C^*$-algebras are then easily
determined from the general theory of graph $C^*$-algebras.
For the quantum real projective space, we also recall the
classification of the classes of irreducible
$*$-representations of its algebra and give a linear basis for this algebra.

\end{abstract}
\newpage\section{Introduction}
In noncommutative geometry, one thinks of quantum spaces as objects dual to
noncommutative algebras in the sense of the Gelfand-Neumark correspondence
between spaces and function algebras.
Analogously to classical topology, where taking quotients of proper group
actions is a standard way of obtaining new
topological spaces, one is then lead to think of Hopf-algebra coinvariant
subalgebras as encoding noncommutative quotient spaces.
For Hopf algebras that are function algebras on finite groups,
their coactions on algebras can be easily understood as group actions
on algebras. Then the coinvariant subalgebras are simply fixed-point
subalgebras.
 In this note, we study two simple examples of
such actions, notably actions of 
  $\Z_2$ on the algebra of the Podle\'s equator quantum sphere $S^2_{q,1}$. 
  For this noncommutative sphere, there is still an analogue of an equator 
formed by the classical points (one-dimensional representations).
One can quotient $S^2_{q,1}$ by the
reflection with respect to the equator plane, and verify that the quotient coincides
with a well-known quantum disc of Klimek and
Lesniewski \cite{kl93}. On the other hand,
one can also define the antipodal action. The quotient under this action is the quantum
real projective
space \cite{h-pm96}. Here we mostly review results from \cite{hms02}
concerning the aforementioned noncommutative quotient spaces, recast in the
more recent language of graph $C^*$-algebras \cite{hs}.

\section{Graph {\boldmath $C^*$}-algebras}
For the convenience of the reader, we briefly recall the definition of a
graph $C^*$-algebra (e.g., see \cite{flr}).
Let $G$ be a countable graph with the set of vertices $G^0$ and 
the set of directed edges
$G^1$.
 (If $e\in G^1$ then $s(e)$ and $r(e)$ are the source and range of 
$e$, respectively.) For the sake of simplicity, we assume that every 
vertex in $G$
emits only finitely many edges. Then $C^*(G)$ is, by definition, the
universal $C^*$-algebra generated by partial isometries
$\{S_e|e\in G^1\}$ with mutually orthogonal ranges and by mutually orthogonal
projections $\{P_v|v\in G^0\}$ such that
\begin{enumerate}
\item
$S_e^*S_e=P_{r(e)}\mbox{ for any }~e\in G^1,$
\item
$P_v=\sum_{s(f)=v}S_fS_f^* ,
\mbox{ for any $v\in G^0$ emitting
at least one edge}$.
\end{enumerate}
These so-called graph $C^*$-algebras generalize the
classical Cuntz-Krieger algebras \cite{ck80}.

According to the general Cuntz formula valid for graph $C^*$-algebras
\cite[Theorem 3.2]{rs}, the $K$-theory of $C^*(G)$ can be calculated
as follows. Let $G^0_+$ be the set of those vertices of a graph $G$ that emit
at least one edge, and let $\Z G^0$ and $\Z G^0_+$ be the free 
abelian groups with generators $G^0$ and $G^0_+$, respectively.
Let $A_G:\Z G^0_+\rightarrow \Z G^0$ be the map defined by
\begin{equation}
\label{ag}
A_G(v)=\left(\mbox{$\sum_{s(e)=v,e\in G^1}$}\;r(e)\right) -v.
\end{equation}
Then
\begin{equation}
\label{k01}
K_0(C^*(G))\simeq {\rm Coker}(A_G),~~~K_1(C^*(G))\simeq {\rm Ker}(A_G).
\end{equation}
\section{Quantum spheres}
The quantum spheres $S^2_{q c}$, $0<|q|<1$, $c\in[0,\infty]$, 
were discovered by Podle\'s as 
$SU_q(2)$-homogeneous spaces \cite{p-p87}.
A uniform description of all these spheres can be given in
the
following way. First, we change the Podle\'s parameter $c\in[0,\infty]$ into
$s\in[0,1]$ via the formula $c=(s^{-1}-s)^{-2}$ (equivalently,
$s=\frac{2\sqrt{c}}{1+\sqrt{1+4c}}$) \cite{bm00}. Then we can rescale the
Podle\'s
generators $A$ and $B$:
\begin{equation}
K:=\left\{\begin{array}{cr}
(1-s^2)A&\mbox{for $s\in[0,1[$ (i.e., $c\in[0,\infty$[)}\\
A&\mbox{for $s=1$ (i.e., $c=\infty$)},\\
\end{array}\right.
\end{equation}
\begin{equation}
L:=\left\{\begin{array}{cr}
(1-s^2)B&\mbox{for $s\in[0,1[$ (i.e., $c\in[0,\infty$[)}\\
B&\mbox{for $s=1$ (i.e., $c=\infty$)}.\\
\end{array}\right.
\end{equation}
This allows us to define the coordinate $*$-algebras of the family of
Podle\'s
quantum spheres $S^2_{q,s}$ as the $*$-algebras  generated by $K$ and $L$ 
satisfying the relations:
\begin{equation}
K=K^*,\;\;\; LK=q^2KL,\;\;\; L^*L+K^2=(1-s^2)K+s^2,\;\;\;
LL^*+q^4K^2=(1-s^2)q^2K+s^2.
\end{equation}
%
The $C^*$-algebra $C(S^2_{q,s})$ can be given as the norm closure of 
$\cO(S^2_{q,s})$.

It seems interesting that the algebra $\cO(S^2_{q,1})$ is isomorphic with
$\cO(SU_{q^2}(2))/\langle b-b^*\rangle$,
where $\cO(SU_{q^2}(2))$ is the polynomial algebra of $SU_{q^2}(2)$
and $b$ is the upper-off-diagonal generator in the fundamental representation
{\scriptsize$\left(\!\!\!\begin{array}{cc}a&b\\-q^{-2}b^*&a^*\end{array}
\!\!\!\right)$}
of $SU_{q^2}(2)$. (Cf.\ the paragraph above Proposition~3.2 in \cite{hs} and 
the relevant considerations in \cite{hl}.)
Indeed, it follows immediately from the defining relations of
$\cO(S^2_{q,1})$ and $\cO(SU_{q^2}(2))$ (e.g., see \cite{w-sl87} and 
put $\nu=q^2,\alpha=a,\gamma=q^{-2}b^*$)
that the assignment
$K\mapsto q^{-2}b$, $L\mapsto a$ defines an algebra epimorphism
$F:\cO(S^2_{q,1})\rightarrow \cO(SU_{q^2}(2))/\langle b-b^*\rangle$.
The injectivity of $F$ follows from the representation theory of these algebras.
The representations $\rho_\pm$ of $\cO(SU_{q^2}(2))$ 
\cite{vs88}
\begin{equation}
\rho_\pm(b)e_k=\pm q^{2(k+1)}e_k,~~~\rho_\pm (a)e_k=\sqrt{1-q^{4k}}\:e_{k-1},~~~
\langle e_k|e_l\rangle=\delta_{kl},~~~ k,l\in\N,
\end{equation}
annihilate the ideal $\langle b-b^*\rangle$ and composed with $F$
agree with Podle\'s representations \cite{p-p87}:
$\pi_\pm=\rho_\pm\circ F$. Since $\pi_+\oplus\pi_-$ is faithful,
one can conclude that $F$ is injective. Hence it is an isomorphism.
The geometric meaning of this isomorphism is that $S^2_{q,1}$
is embedded as an equator in $SU_{q^2}(2)$ thought of as a
quantum 3-sphere. The other extreme value of $s$, i.e., $s=0$,
also simplifies the relations. In this case, it turns out that
$\O(S^2_{q,0})$ is isomorphic with the fixed-point subalgebra
$\O(SU_{q}(2))^{U(1)}$, so that we can interpret $S^2_{q,0}$ as the
quotient sphere $SU_{q}(2)/{U(1)}$ in the spirit of the Hopf fibration.
(Here the action of $U(1)$ is given by rescaling the generator $a$
by $e^{i\phi}$ and $b$ by $-e^{i\phi}$.)
This way we can view the family of Podle\'s spheres $S^2_{q,s}$
as an approximation between the quotient sphere
$S^2_{q,0}\simeq SU_{q}(2)/{U(1)}$ and the embedded sphere
$S^2_{q,1}\subset SU_{q^2}(2)$. Although the desired $\Z_2$-actions
can be defined on the $C^*$-level for any $s\in ]0,1]$
(for $s>0$, $C(S^2_{q,s})\simeq C(S^2_{q,1})$),
in what follows, we restrict our attention to $S^2_{q,1}$
because this is where we can define these actions on the algebraic level.

It is shown in \cite[Proposition 3.1]{hs} that, for any $q\in]0,1[$, 
there exists a $C^*$-algebra isomorphism $\phi_q$ from $C(S^2_{q,1})$
to the graph $C^*$-algebra $C^*(G_1)$ corresponding to 
the following graph: 
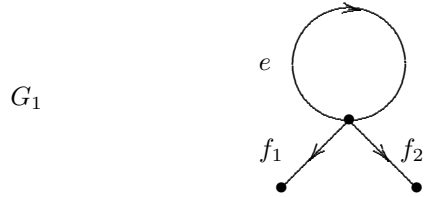
\[ \beginpicture
\setcoordinatesystem units <1.5cm,1.5cm>
\setplotarea x from -5 to 5, y from -1 to 1.1  
\put {$\bullet$} at -0.6 -0.6
\put {$\bullet$} at  0 0 
\put {$\bullet$} at 0.6 -0.6 
\circulararc 360 degrees from 0 0 center at 0 0.5 
\setlinear 
\plot 0 0 -0.6 -0.6 / 
\plot 0 0 0.6 -0.6 /
\put {$e$} [l] at -0.8 0.5
\put {$f_1$} [l] at -0.8 -0.27
\put {$f_2$} [l] at 0.45 -0.27
\put {$G_1$} [l] at -3 0.2
\arrow <0.235cm> [0.2,0.6] from 0 0.99 to 0.1 0.975 
\arrow <0.235cm> [0.2,0.6] from -0.345 -0.351 to -0.35 -0.357 
\arrow <0.235cm> [0.2,0.6] from 0.345 -0.354 to 0.35 -0.36 
\endpicture \]
(Note that $C^*(G_1)$ is generated by partial isometries $S_e$, $S_{f_1}$ and 
$S_{f_2}$, corresponding to the edges of $G_1$.)
The matrix of the group homomorphism $A_{G_1}$ (see (\ref{ag}))
is 
\begin{equation}
A_{G_1}=\left(\begin{array}{c}0\\1\\1\end{array}\right):\Z\longrightarrow\Z^3.
\end{equation}
Hence the formulas (\ref{k01}) yield $K_0(C(S^2_{q,1}))\simeq \Z^2$,
$K_1(C(S^2_{q,1}))\simeq 0$, which agrees with the computation of
\cite{mnw91} and with the classical situation.

Note that the primitive ideal space $\Prim(C(S^2_{q,1}))$ (coinciding with
$\Prim(C(S^2_{q,s}))$ for all $s\in ]0,1]$) can be described as follows: One has a
circle $S^1$ with the usual topology and two extra points, which are separated
from each other but can not be separated from the circle. This is easily
derived from the irreducible representations \cite{p-p87} or using general
arguments of the theory of graph $C^*$-algebras.

\section{Quantum disc} 
The $*$-algebra of the quantum disc $D_q$ is defined as follows:
\begin{equation}
{\cal O}(D_q):=\C\langle x,x^*\rangle/\langle x^*x-q xx^*-(1-q)\rangle,~~0<q<1.
\end{equation}
This is a one-parameter subfamily of a two-parameter family of quantum discs
introduced by Klimek and Lesniewski \cite{kl93} as homogeneous spaces of
$SU_q(1,1)$. It was shown in \cite{kl93} that
$\|\pi(x)\|=1$ in every bounded $*$-representation of ${\cal O}(D_q)$. 
Thus the $C^*$-algebra $C(D_q)$ can be defined using the supremum
over the norms in bounded $*$-representations. (Note here that there are
unbounded representations -- the disc algebra is easily transformed into the
$q$-oscillator algebra.) As also shown in \cite{kl93},
$C(D_q)$ is isomorphic to the $C^*$-algebra of the one-sided shift
(Toeplitz algebra). 
This isomorphism is provided by a faithful infinite dimensional 
representation 
$\pi$ \cite[p.14]{kl93}.
The one-dimensional
representations of $C(D_q)$ form a circle that can be considered the
boundary of the quantum disc $D_q$. The infinite dimensional
representation corresponds to the interior of that disc.

We define now a $\Z_2$-action by sending $1\in \Z_2$ to 
$r_1: C(S^2_{q,1})\rightarrow C(S^2_{q,1})$ that 
``identifies
upper and lower hemispheres'': 
\begin{equation}
r_1(L)=L,~~r_1(K)=-K.
\end{equation}
We have shown in \cite{hms02} the following results concerning this action:
\begin{itemize}
\item
The fixed-point polynomial algebra
${\cal O}(S^2_{q,1}/\Z_2):=\{a\in {\cal O}(S^2_{q,1})~|~r_1(a)=a\}$
is the $*-$subalgebra
generated by $L$, and can be identified with ${\cal O}(D_{q^4})$ 
by sending the generator $x$ to the generator $L$.
This extends to the polynomial level the fact that the $r_1$-fixed
point subalgebra of $C(S^2_{q,1})$ coincides with the Toeplitz algebra.
\item
The $\Z_2$-extension ${\cal O}(D_{q^4})\subset {\O(S^2_{q,1})}$ defined in
this way is not Galois. (The $\Z_2$-action is not free.)
\item
All $*$-representations of ${\cal O}(D_{q^4})$
are restrictions of $*$-representations of ${\cal O}(S^2_{q,1})$.
\item
The automorphism $r_1$ commutes with the $SU_q(2)$-induced
$U(1)$-action on $S^2_{q,1}$.
\end{itemize}

Under the aforementioned isomorphism 
$\phi_q:C(S^2_{q,1})
\rightarrow C^*(G_1)$, the order two automorphism 
$r_1$ of $C(S^2_{q,1})$ is 
transformed into the automorphism of $C^*(G_1)$ determined by 
\begin{equation}
 S_e\mapsto S_e, \;\;\; S_{f_i}\mapsto S_{f_{3-i}}. 
 \end{equation}
This automorphism of $C^*(G_1)$ is induced from the automorphism of the graph
$G_1$ that fixes the edge $e$ and interchanges $f_1$ with $f_2$.
The fixed-point subalgebra for this $\Z_2$-action coincides with
 the Toeplitz algebra and 
corresponds to the following graph:  
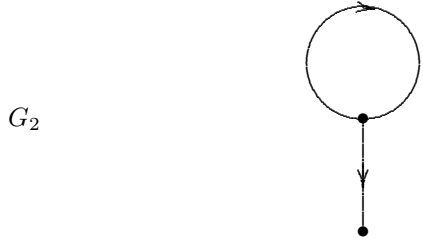
\[ \beginpicture
\setcoordinatesystem units <1.5cm,1.5cm>
\setplotarea x from -5 to 5, y from -1 to 1.1  
\put {$G_2$} at -3 0
\put {$\bullet$} at  0 0 
\put {$\bullet$} at 0 -1 
\circulararc 360 degrees from 0 0 center at 0 0.5 
\setlinear 
\plot 0 0  0 -1 / 
\arrow <0.235cm> [0.2,0.6] from 0 0.99 to 0.1 0.975 
\arrow <0.235cm> [0.2,0.6] from 0 -0.4 to 0 -0.6 
\endpicture \] 
The $K$-theory of the Toeplitz algebra is well-known (e.g. see \cite{w-ne93}).
Here, we can directly determine it from (\ref{ag}) and (\ref{k01}):
\begin{equation}
A_{G_2}=\left(\begin{array}{c}0\\1\end{array}\right):\Z\longrightarrow\Z^2,~~~
K_0(C(D_q))\simeq\Z,~~~ K_1(C(D_q))\simeq 0.
\end{equation}

\section{Quantum two-dimensional real projective space}
We define the antipodal $\Z_2$-action on $S^2_{q,1}$ by
sending $1\in\Z_2$ to the $*$-automorphism $r_2$ of 
$C(S^2_{q,1})$ defined by
\[
r_2(K)=-K,~~r_2(L)=-L.
\]
The $*$-algebra ${\cal O}(\R P^2_q):=\{a\in {\cal
O}(S^2_{q,1})|r_2(a)=a\}$ 
describes the quantum two-dimensional real projective 
space $\R P^2_q$.
This algebra is generated by
$
P=K^2,~~R=L^2,~~T=KL.
$
They fulfill the relations
\begin{equation}\label{p}
P=P^*,
\end{equation}
\begin{equation}\label{pc}
RP=q^8PR,~~RT=q^4TR,~~PT=q^{-4}TP,
\end{equation}
\begin{equation}\label{pp}
T^2=q^2PR,~~RT^*=q^2T(-q^4P+I),~~R^*T=q^{-2}T^*(-P+I),
\end{equation}
\begin{equation}\label{pcpr}
RR^*=q^{12}P^2-q^4(1+q^4)P+I,~~R^*R=q^{-4}P^2-(1+q^{-4})P+I,
\end{equation}
\begin{equation}\label{pcpt}
TT^*=-q^4P^2+P,~~T^*T=q^{-4}(P-P^2).
\end{equation}
One can show that
${\cal O}(\R P^2_q)$ is isomorphic to the free $*$-algebra 
generated by $P,~R,~T$ divided by the ideal of relations 
(\ref{p})--(\ref{pcpt}), see \cite{hms02}.
Concerning representations, we have:
\begin{theorem}{Theorem}$\!\!\!$\cite{hms02}$\;$
There are no unbounded $*$-representations of the $*$-algebra
${\cal O}(\R P^2_q)$. Up to unitary equivalence, all irreducible
$*$-representations of this algebra  are the following:

(i) A family of one-dimensional representations
$\rho_\theta:{\cal O}(\R P^2_q)\rightarrow \C$, parameterised by 
$\theta\in[0,2\pi)$,
given by
\begin{equation}
\rho_\theta(P)=\rho_\theta(T)=0,~~\rho_\theta(R)=e^{i\theta}.
\end{equation}

(ii) An infinite dimensional representation $\rho$ (in a Hilbert
space $H$ with an orthonormal basis $\{e_k\}_{k\in\N}$) given by
\begin{eqnarray}
\label{rp}\rho(P)e_k&=&q^{4k}e_k,\\
\label{rt}\rho(T)e_k&=&q^{2(k-1)}\sqrt{1-q^{4k}}e_{k-1},\\
\label{rr}\rho(R)e_k&=&\sqrt{(1-q^{4k})(1-q^{4(k-1)})}e_{k-2}.
\end{eqnarray}
\end{theorem}
{\sc Proof (sketch):}
The boundedness of $\rho(P)=\rho(P^*)$ comes from the relation
$T^*T=q^{-4}(P-P^2)$. Then it is an
easy consequence of the relations (\ref{pcpr}) and (\ref{pcpt}) that
$\rho(R)$ and $\rho(T)$ are also bounded.
Next, the assumption $\rho(P)=0$ immediately leads to $\rho(T)=0$,
and the only remaining relations are $\rho(R^*R)=1=\rho(RR^*)$.
This yields the one-dimensional representations.
On the other hand,
assuming $\rho(P)\neq 0$ implies, by the irreducibility of $\rho$, that
${\rm Ker}\rho(P)=0$. Using the characterisation of the spectrum by
approximate eigenvectors, and taking advantage of the relations,
it is possible to identify the 
spectrum of $\rho(P)$ with $\{0\}\cup\{q^{4k}|k\in\N\}$.
 Now, one builds a Hilbert space out of
the eigenvectors of $\rho(P)$, and identifies $\rho(T)$ and $\rho(R)$
as weighted-shift operators. \hfill$\Box$

Note that  the irreducible $*$-representations of ${\cal O}(\R P^2_q)$ are the
restrictions of the irreducible $*$-representations
of $\O(S^2_{q,1})$. 
\bpr
(i) The infinite dimensional representation $\rho$ of $\O(\R P^2_q)$
is faithful.\\
(ii) The set
$\{P^kR^l~|~ k,l\in\N\}\cup \{P^k{R^*}^l~|~ k\in \N, l\in\N\setminus\{0\}\}
\cup \{P^kR^lT~|~ k,l\in\N\}\cup\{P^k{R^*}^lT^*~|~ k,l\in\N\}$
is a basis of the
vector space $\O(\R P^2_q)$.
\epr
{\sc Proof:} 
We know from the proof of \cite[Proposition 4.3]{hms02} that the set (ii)
generates $\O(\R P^2_q)$ as a vector space. Thus, a general element of 
$\O(\R P^2_q)$ is of the form

\begin{equation}
x=\sum_{k,l\geq 0}a_{kl}P^kR^lT+\sum_{k,l\geq 0}b_{kl}P^k{R^*}^lT^*+
\sum_{k,l\geq 0}c_{kl}P^kR^l+\sum_{k\geq 0,l\geq 1}d_{kl}P^k{R^*}^l.
\end{equation}
For the proof of both claims of the proposition it is sufficient to show
 that it follows 
from $\rho(x)=0$ that
all the coefficients $a_{kl}, b_{kl}, c_{kl}, d_{kl}$ vanish.
Acting with $\rho(x)$ onto a basis vector $e_n$ of the representation space,
making use of (18) -- (20) and the formulas 
\begin{eqnarray*}
\rho(T^*)e_n&=&q^{2n}\sqrt{1-q^{4(n+1)}}e_{n+1},\\
\rho(R^*)e_n&=&\sqrt{(1-q^{4(n+1)})(1-q^{4(n+2)})}e_{n+2},
\end{eqnarray*}
we obtain
\begin{eqnarray*}
\rho(x)e_n&=&\sum_{k\geq 0,l\leq\frac{n-1}{2}}a_{kl}q^{4k(n-1-2l)}
\sqrt{(1-q^{4(n-2l)})\cdots
(1-q^{4(n-1)})}q^{2(n-1)}\sqrt{1-q^{4n}}e_{n-1-2l}\\
&&+\sum_{k,l\geq 0}b_{kl}q^{4k(n+1+2l)}\sqrt{(1-q^{4(n+1+2l)})
\cdots(1-q^{4(n+2)})}q^{2n}\sqrt{1-q^{4(n+1)}}e_{n+1+2l}\\
&&+\sum_{k\geq 0,l\leq\frac{n}{2}}c_{kl}q^{4k(n-2l)}\sqrt{(1-q^{4(n-2l+1)})
\cdots(1-q^{4n})}e_{n-2l}\\
&&+\sum_{k\geq 0,l\geq 1}d_{kl}q^{4k(n+2l)}
\sqrt{(1-q^{4(n+2l)})\cdots(1-q^{4(n+1)})}e_{n+2l}\\
&=&0.
\end{eqnarray*}
Fixing $l$ we get the following four sets of equations:
\begin{eqnarray}\label{a}
\sum_{k\geq 0}a_{kl}q^{4k(n-1-2l)}&=&0,~~~n\geq 1,~l\leq\frac{n-1}{2},
\\
\label{b}
\sum_{k\geq 0}b_{kl}q^{4k(n+1+2l)}&=&0,~~~n\geq 0,~l\geq 0,
\\
\label{c}
\sum_{k\geq 0}c_{kl}q^{4k(n-2l)}&=&0,~~~n\geq 0,~l\leq\frac{n}{2},
\\
\label{d}
\sum_{k\geq 0}d_{kl}q^{4k(n+2l)}&=&0,~~~n\geq 0,~l\geq 1.
\end{eqnarray}
Let us consider the equations (\ref{b}) for fixed $l\geq 0$. Since $b_{kl}$ is
different from 0 only for a finite number of indices $k$, there is a biggest
$k=N$ with $b_{kl}\neq 0$. Then the first $N$ (i.e., $n=0,\ldots, N$)
 of the equations (\ref{b})
are a system of linear equations for $b_{0l},\ldots,b_{Nl}$ whose matrix
is a Vandermonde matrix
\[
\left(\begin{array}{cllll}
1&x_0&x_0^2&\ldots&x_0^N\\
1&x_1&x_1^2&\ldots&x_1^N\\
\vdots&\vdots&\vdots&&\vdots\\
1&x_N&x_N^2&\ldots&x_N^N
\end{array}\right)
\]
with $x_i=q^{4(i+1+2l)}$. Since $x_i\neq x_j$ for $i\neq j$, the
determinant $\prod_{0\leq i<j\leq N}(x_j-x_i)$ of this matrix is nonzero.
It follows that $b_{kl}=0$. One argues analogously for the remaining
equations (\ref{a}), (\ref{c}) and (\ref{d}). Note that for (\ref{a})
and (\ref{c}) one has to consider a set of equations for
$n=n_0,\ldots,n_0+N$ with $n_0$ such that $l\leq\frac{n_0-1}{2}$ and
$l\leq\frac{n_0}{2}$ respectively. Again one arrives at Vandermonde
matrices with $x_i=q^{4(n_0+i-1-2l)},q^{4(n_0+i+2l)},q^{4(i+2l)}$
for the equations (\ref{a},\ref{c},\ref{d}) respectively. The determinants
of these matrices are nonzero, which leads to the desired conclusion.
\epf 
 
It follows that the $*$-algebra $\O(\R P_q^2)$ is faithfully embedded in
its enveloping $C^*$-algebra $C(\R P_q^2)$. 
Considering special simple elements of the algebra, it is easily seen
that the kernels of the representations $\rho_\theta$ are not contained in
each other for different $\theta$. On the other hand, the kernel of $\rho$
is trivial and thus contained in each $\ker\rho_\theta$. Therefore, the
primitive ideal space $\Prim(C(\R P_q^2))$ consists of a circle $S^1$,
which has the usual topology, and an extra point that can not be separated
from that circle. This picture also follows from general results of the
theory of graph $C^*$-algebras. Note that $\Prim(C(D_q))=\Prim(C(\R P_q^2))$,
the different topologies of the two quantum spaces (different $K$-groups) are
not visible in the primitive ideal space.

Next, let us observe that,
since $r_2$ and the right coaction of $SU_q(2)$ on $S^2_{q,1}$ 
commute \cite{p-p87,hms02}, $\R P^2_q$ is a quantum homogeneous 
space of $SU_q(2)$.
One can also show that
the $\Z_2$-extension ${\cal O}(\R P^2_q)\subset {\cal O}(S^2_{q,1})$
is Galois \cite{h-pm96}, i.e., we have an algebraic quantum principal bundle.
Moreover, since $\O(S^2_{q,1})\otimes {\rm Map}(\Z_2,\C)$
is norm dense in $C(S^2_{q,1})\otimes {\rm Map}(\Z_2,\C)$,
one can see that the antipodal $\Z_2$-action on $C(S^2_{q,1})$
is principal in the sense of \cite{e-da00}.

Finally, let us look at $C(\R P^2_q)$ as a graph $C^*$-algebra.
The
automorphism $r_2$ of $C(S^2_{q,1})$  corresponds via the 
previously-used isomorphism
$\phi_q$ to the $\Z_2$-action  on $C^*(G_1)$ determined by 
$$ S_e\mapsto -S_e, \;\;\; S_{f_i}\mapsto -S_{f_{3-i}}. $$ 
This is a quasi-free automorphism of $C^*(G_1)$ not induced from a graph
automorphism.
The fixed-point subalgebra for this action coincides with $C(\R P^2_q)$ and
is isomorphic with the $C^*$-algebra of the following graph \cite{hs}: 
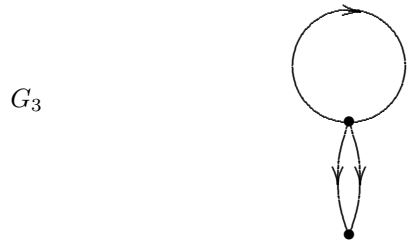
\[ \beginpicture
\setcoordinatesystem units <1.5cm,1.5cm>
\setplotarea x from -5 to 5, y from -0.5 to 1.1 
\put {$\bullet$} at  0 -1
\put {$\bullet$} at  0 0 
\circulararc 360 degrees from 0 0 center at 0 0.5 
\setquadratic 
\plot 0 0 -0.1 -0.5 0 -1 / 
\plot 0 0 0.1 -0.5 0 -1  /
\put {$G_3$} [l] at -3 0.2
\arrow <0.235cm> [0.2,0.6] from 0 0.99 to 0.1 0.983
\arrow <0.235cm> [0.2,0.6] from -0.1 -0.54 to -0.1 -0.55 
\arrow <0.235cm> [0.2,0.6] from 0.1 -0.54 to 0.1 -0.55 
\endpicture \] 
The group homomorphism $A_{G_3}$ (see (\ref{ag}))
has the form
\begin{equation}
A_{G_3}=\left(\begin{array}{c}0\\2\end{array}\right):\Z\longrightarrow\Z^2.
\end{equation}
It follows now from (\ref{k01}) that the $K$-groups of $C(\R P^2_q)$
agree with their classical counterparts:
\begin{theorem}{Theorem}$\!\!\!$\cite{hms02}$\;$
$
K_0(C(\R P^2_q))\simeq \Z\oplus\Z_2,~~ K_1(C(\R P^2_q))\simeq 0.
$
\end{theorem}

\noindent
{\bf Acknowledgements:}
The work on this paper was partially supported by a Marie Curie fellowship
(PMH),  Deutsche Forschungsgemeinschaft (RM),
the Research Management Committee of the University of Newcastle and
Max-Planck-Institut f\"ur Mathematik Bonn (WS).
All authors are grateful to Mathematisches Forschungsinstitut Oberwolfach for support via its
Research in Pairs programme. It is a pleasure to thank Giovanni Landi for a very
helpful discussion.

\end{document}